\documentclass[12pt]{article}
\usepackage{amssymb}
\usepackage{enumerate}
\usepackage{graphicx}
\usepackage{amsthm}
\usepackage{color}

\def\R{{\rm I\! R}}

\def\dt{{ d \over dt}}

\newtheorem{theorem}{Theorem}

\title{An extension to the planar Markus-Yamabe Jacobian conjecture }

\author{Marco Sabatini 
\footnote{Dipartimento di Matematica, Univ. di Trento, I-38123 Povo (TN) - Italy; email: marco.sabatini@unitn.it. This paper has been partially supported by GNAMPA, Gruppo Nazionale per l'Analisi Matematica, la Probabilit\`a e le loro Applicazioni. }
}
\begin{document}
\maketitle
\begin{abstract}  We extend the planar Markus-Yamabe Jacobian Conjecture to differential systems having jacobian matrix with eigenvalues with negative or zero real parts.

 \end{abstract}

\section{Introduction}

Let
\begin{equation}  \label{sysn}
\dot X = F(X), \qquad X \in \R^n, \qquad F\in C^1(\R^n,\R^n)
\end{equation}
be a first order differential system. Let us denote by $J_F(X)$ the jacobian matrix of $F(X)$. If $O$ is a critical point of (\ref{sysn})   and  the eigenvalues of $J_F(O)$ have negative real parts, then $O$ is asymptotically stable \cite{BS}. In particular, all orbits starting close enough to $O$ tend asymptotically to $O$.

In \cite{MY} the question was raised, whether $J_F(X)$ having eigenvalues with negative real parts for every $X \in \R^n$ imply $O$ to be globally asymptotically stable, i. e. whether all orbits in $\R^2$ tend asymptotically to $O$.  Such a problem was named {\it Markus-Yamabe Jacobian Conjecture} and several results were obtained under various additional hypotheses. A key step was made in \cite{Ol}, where it was proved that under Markus-Yamabe hypotheses, for planar systems the global asymptotic stability of $O$ is equivalent to the injectivity of $F(X)$. Such a result led to study the problem applying methods previously used to study injectivity. 
 The Markus-Yamabe Jacobian Conjecture was solved in the positive in \cite{Fe,Gl,Gu} for planar systems, and was proved to have negative answer in higher dimensions \cite{BL,CEGHM}. 
The three approaches proposed in in \cite{Fe,Gl,Gu}  first prove the injectivity  of $F(X)$, then as a consequence get the global asymptotic stability. Actually, in all such papers injectivity is proved under much weaker hypotheses than that of negative real parts. In fact, it is sufficient to assume that the Jacobian matrix has nowhere real positive eigenvalues. 

Such general results did not lead to similarly general results in the study of the systems dynamics. This is likely due to the fact that accepting the possibility of eigenvalues with different real parts (positive, zero or negative) at different points of the plane does not allow to apply the procedure developped in \cite{Ol} to establish the equivalence of injectiviy and global asymptotic stability. On the other hand, eigenvalues with zero real parts are compatible with asymptotic stability, even if not sufficient to imply it.

In this paper we  assume $J_F(X)$ to be non-singular and have eigenvalues with non-positive real parts for all $X \in \R^2$. Differently from the classical case, in this case a system does not necessarily have a  globally asymptotically stable critical point. If a critical point exists, we prove that either such  a system has a global center, or there exists a globally asymptotically stable compact set. We show by an example that such a global attractor is not necessarily a critical point.   If the system is analytic the conclusion can be sharpened, proving that either there exists a global center, or a globally asymptotically stable critical point.  Our results follow from Olech approach to global attractivity \cite{Ol} and Fessler theorem about global injectivity \cite{Fe}.

\section{Results}

We consider  maps $F \in C^1(\R^2 , \R^2)$, $F(x,y) = (P(x,y),Q(x,y))$. We denote partial derivatives by subscripts. Let
$$
J_F (x,y) = \left( \matrix{P_x(x,y) & P_y(x,y) \cr Q_x(x,y) & Q_y(x,y)}  \right)
$$ 
be tha jacobian matrix of $F$ at $(x,y)$. We denote by $D(x,y) = \det J_F(x,y)  = P_x(x,y) Q_y(x,y) - P_y(x,y) Q_x(x,y)$ its determinant and by $T(x,y) = P_x(x,y) + Q_y(x,y)$ its trace. $T(x,y)$ is the divergence of the vector field $F(x,y)$.

In what follows we consider the differential system associated to $F$:
\begin{equation}  \label{sys1} 
 \left\{ \matrix{\dot x =  P(x,y) \cr \dot y = Q(x,y) .}  \right. 
\end{equation}
We denote by $\phi(t,x,y) $ the local flow defined by (\ref{sys1}).   
We say that a critical point $O$ of (\ref{sys1}) is a {\it center} if it has a punctured neighbourhood filled with non-trivial cycles surrounding $O$. The largest connected set $N_O$ filled with such cycles is called {\it period annulus of $O$}. If $N_O = \R^2 \setminus \{ O \}$, then $O$ is said to be a {\it global center}.  
We say that a critical point $O$ of (\ref{sys1}) is {\it asymptotically stable} if it is stable and attractive \cite{BS}. In this case we denote by  $A_O$ its {\it attraction region}. If $A_O = \R^2$ then $O$ is said to be {\it globally asymptotically stable}.

In the proof of theorem \ref{teorema} we repeatedly use $F$ injectivity. We report here the theorem applied, proved in \cite{Fe}.

\begin{theorem}\label{fessler} 
Let $F \in C^1(\R^2,\R^2)$ be such that:
\begin{itemize}
\item[1)] $D(x,y) > 0$ for all $(x,y) \in \R^2$;
\item[2)] there is a compact set $K \subset \R^2$ such that $J_F (x,y)$ has no real positive eigenvalues for any $(x,y) \not\in K$.
\end{itemize}
Then $F$ is injective.
\end{theorem}

For the sake of simplicity, without loss of generality from now on we assume $O=(0,0)$. The  hypotheses we consider rely only on  derivatives properties, hence they do not change after a translation. We set
$$
T_- = \{  (x,y) : T(x,y) < 0 \},
$$
and denote by $\overline{T_-}$ its closure. We denote by   $\mu$  the 2-dimensional Lebesgue measure.

\begin{theorem}\label{teorema} 
Assume $D(x,y) > 0$ and $T(x,y) \leq 0$ for all $(x,y) \in \R^2$. Let $O$ be a critical point of  (\ref{sys1}). Then: 
\begin{itemize}
\item[i)]  $O$ is a center if and only if  it has a neighbourhood $U_O$ such that $T(x,y)$ vanishes identically on $U_O$; in such a  case (\ref{sys1}) is Hamiltonian on all of $N_O$; if,  additionally,  $F$ is analytic, then the system is Hamiltonian and $O$ is a global center.\item[ii)]  
$O$ is asymptotically stable if and only if it belongs to  $\overline{T_-}$;  in such a  case  $O$ is globally asymptotically stable.
\item[iii)]  If  $T(x,y) $ does not vanish identically, then there exists a globally asymptotically stable compact set $M$.   

\end{itemize}
\end{theorem}
{\it Proof.}   

$i.1)$
We claim that if $O$ is a center, then $T(x,y) $ vanishes identically on $N_O$. By absurd, assume $T(x^*,y^*) < 0$ for some $(x^*,y^*) \in N_O$. By continuity there exists a neighbourhood $U^*$ of $(x^*,y^*)$ such that $T(x,y) < 0$ for all $(x,y) \in U^*$. Let $\gamma^*$ be the cycle passing through $(x^*,y^*)$ and $\Delta^*$ the bounded planar region having $\gamma^*$ as boundary.  $\Delta^*$ is invariant, hence $\mu(\Delta^*) = \mu (\phi(t,\Delta^*))$ for all $t\in \R$. By Liouville theorem one has
$$
0 = \dt\, \mu (\phi(t,\Delta^*)) = \int_{\phi(t,\Delta^*)} T(x,y)\, dx\, dy < 0,
$$
because $T(x,y) < 0 $ on $\phi \left(t,\Delta^* \cap  U^*\right)$, contradiction. 

$i.2)$
Vice-versa, assume $T(x,y) $ to vanish identically on a neighbourhood $U_O$ of $O$. Then the system is Hamiltonian on a simply connected neighbourhood $V_O \subset U_O$.  Let $H(x,y)$ be its Hamiltonian function. One has
\begin{equation}  \label{sysham} 
 \left\{ \matrix{\dot x =  P(x,y) = -H_y(x,y)\cr \dot y = Q(x,y) =  H_x(x,y) \hfill.}  \right. 
\end{equation}
The Hessian matrix of $H(x,y)$ is
$$
\left( \matrix{H_{xx} & H_{xy} \cr H_{yx} & H_{yy}}  \right) = \left( \matrix{Q_x& Q_y \cr -P_x & -P_y }  \right) .
$$ 
The Hessian determinant is $H_{xx} H_{yy} - H_{xy} H_{yx}= P_x Q_y - P_y Q_x = D(x,y) > 0$, hence $H(x,y)$ has a minimum at $O$. As a consequence, $O$ is a center. 

$i.3)$
If additionally $F$ is analytic, then also $T(x,y)$ is analytic. If it vanishes in a neighbourhood of $O$ then it vanishes on all of $\R^2$, hence the system is Hamiltonian on all of $\R^2$. 
We claim that $N_O$ is unbounded. In fact, let us assume by absurd $N_O$ is bounded, hence also $\partial N_O$ is bounded. By $F$'s injectivity \cite{Fe}, $\partial N_O$ contains no critical points, hence by Poincar\'e-Bendixson theorem $\partial N_O$ is a non-trivial cycle.  One can consider the Poincar\'e map defined on a section $\Sigma$ of $\partial N_O$. Such a map is analytic and coincides with the identity map on $\Sigma\, \cap\, N_O $ , hence it coincides with the identity map on all of $\Sigma$. As a consequence every orbit meeting $\Sigma\, \cap\, \partial N_O$ is a cycle, hence $\partial N_O $ is contained in the period annulus, contradicting the fact that it is the boundary of $N_O$. 
Moreover, every connected components of $\partial N_O $ is unbounded. In fact, if a connected components of $\partial N_O $ was bounded, then by its invariance and by Poincar\'e-Bendixson theorem either it would be a cycle or it would contain a critical point. The former case has already been considered above, the latter one can be excluded by the injectivity of $F$.

$i.4)$ In order to prove that $O$ is a global center we use again the injectivity of $F$. For  $\varepsilon > 0$ let $B_\varepsilon$ be the open disk of radius $\varepsilon > 0$ centered at $O$. $F$ is a diffeomorphism, hence the anti-image $ D_\varepsilon = F^{-1}\left(  B_\varepsilon \right)$ is an open neighbourhood of $O$. By construction  and by the injectivity of $F$, $D_\varepsilon$ contains all the points of $(x,y) \in \R^2$ such that $|F(x,y) | <  \varepsilon $, hence for all $(x,y) \not\in  D_\varepsilon $ one has $|F(x,y) | \geq  \varepsilon $.  Let us choose $\varepsilon $ small enough such that  $ \partial N_O \cap  D_\varepsilon = \emptyset$. Let $\partial N_O^u$ be an unbounded component of $\partial N_O$. 
Then working as in \cite{Ol}, since $T(x,y) \leq 0$ and $|F(x,y) | \geq  \varepsilon $ outside  $ D_\varepsilon$, one proves that  every orbit starting close enough to $\partial N_O^u$ is unbounded too, hence it is not a cycle, contradicting the fact that $\partial N_O^u$ is in the boundary of $N_O$. As a consequence $\partial N_O = \emptyset$ and $N_O = \R^2 \setminus \{ O \}$. 

$ii)$ 
Assume $O$ to be asymptotically stable and $A_O$ its region of attraction. By hypothesis, in every neighbourhood of $O$ there are points such that $T(x,y) <0$, and by continuity this occurs in an open subset of $A_O$.  If by absurd $A_O$ is bounded, then by its invariance, for all $t$
$$
0 = \dt\, \mu (\phi(t,A_O)) = \int_{\phi(t,A_O)} T(x,y)\, dx\, dy < 0,
$$
 contradiction. Hence   $A_O$ is unbounded. Assume by absurd there exists a bounded connected component $\partial A_O^b$ of $\partial A_O$. As above, by Poincar\'e-Bendixson theorem either it is a cycle or contains  a critical point. If it is a cycle, it cannot surround $O$, since in such a case $A_O$ would be bounded. Hence  it  surrounds another critical point, violating $F$ injectivity. The same violation would occur if $\partial A_O^b$ contained a critical point. Then the argument proceeds as in point $i.4)$, showing that $\partial A_O = \emptyset$ and $A_O = \R^2$.

Vice-versa, assume  $O \in\overline{T_-}$. Then $T(x,y) $ does not vanish identically on any neighbourhood $U_O$ of $O$, hence by point $i)$ it is not a center. By the hypotheses on $D(0,0)$ and $T(0,0)$, $O$ is a non degenerate elementary critical point of center-focus type, according to the real part of its eigenvalues. If such real parts are negative $O$ is a focus, hence asympotically stable. If such real parts are zero, one proves, as at the beginning of point $i)$, that $O$ cannot be accumulation point of cycles, hence it is  asympotically stable. Working as in point $i.4)$ one  proves that it is globally asymptotically stable.

$iii)$  If $O \in \overline{T_-}$, then point $ii)$ applies and one can take $M = \{ O \}$. \\
If $O \not\in\overline{T_-}$, it has a neighbourhood $U_O$ where $T(x,y)$ vanishes identically, hence it is a center. We claim that $N_O$ is bounded. In fact, if $N_O$ is unbounded one can proceed as in point $i.4)$, in order to prove that every orbit starting close enough to $\partial N_O^u$ is unbounded, contradicting the fact that $\partial N_O^u$ is part of the boundary. The boundedness of $N_O$ implies the boundeness of $\partial N_O$, which is a cycle, by the absence of critical points on $\partial N_O^u$.  Let us consider a section $\Sigma$ of $\partial N_O$ and its Poincar\'e map. Such a map is the identity on $\Sigma \cap N_O$, and has no fixed points on $\Sigma \setminus N_O$, otherwise there would be a cycle $\gamma$ containing $\partial N_O$, $T(x,y)$ would vanish identically inside $\gamma$  and every orbit inside $\gamma$ would be a cycle, contradicting the fact that $\partial N_O$ is the boundary of $N_O$. Hence the Poincar\'e map is strictly monotone, which implies either attractivity or repulsivity of $\partial N_O$. Repulsivity is not compatible with the sign of the divergence, hence $\partial N_O$ is attractive, and $\overline{N_O}$ is asymptotically stable. Its global attractivity can be proved as in $i.4)$ and $ii)$, proving that the boundary of its region of attraction is empty.

 \hfill  $\clubsuit$

An example of globally asymptotically stable critical point belonging to $ \overline{T_-}$
is the origin in the following differential system,
\begin{equation}  \label{sys2} 
 \left\{ \matrix{\dot x =  y   \hfill \cr \dot y = - x - y^3 }  \right. 
\end{equation}
for which one has
$$
J_F (x,y) = \left( \matrix{0 & 1 \cr -1 & -3y^2}  \right) .
$$ 
One has $D(x,y) = 1$, $T(x,y) = - 3 y^2 \leq 0$, hence $T_-$ is $x$-axis.

If (\ref{sys1}) is not analytic, then a center need not be global. We construct a system satisfying the hypotheses of theorem \ref{teorema}, having a non-global center and a globally asymptotically stable compact set.
Let $\alpha \in C^\infty(\R,\R)$ be such that 
$$
 \left\{ \matrix{\alpha(r) = 0,  & r \leq 1 \cr \alpha(r) > 0,  & r > 1  \cr \alpha'(r) > 0, & r >1 .}  \right. 
$$
Let us set $r = \sqrt{x^2 + y^2}$. The vector field defined by the system
\begin{equation}  \label{sysr} 
 \left\{ \matrix{\dot x =  y - x\, \alpha(r)  \hfill \cr \dot y = - x - y\, \alpha(r) }  \right. 
\end{equation}
Setting $\displaystyle{cr = x, sr = y}$,   the jacobian matrix of the vector field is
$$
J (x,y) = \left( \matrix{-\alpha(r) - xc\alpha'(r)  & 1 - xs \alpha'(r) \cr -1 - yc \alpha'(r) & -\alpha(r) - ys \alpha'(r) }  \right) .
$$ 
Its determinant is $1 + \alpha^2(r) + r \alpha(r)\alpha'(r) > 0$ and its trace is $-2 \alpha(r) - 2 r \alpha'(r) \leq 0$. For $r \leq 1$ the trace is zero, for $r > 1$ the trace is negative.  The system (\ref{sysr}) is Hamiltonian for $r \leq 1$, with a center at $O$ whose central region is the disk of radius 1 centered at $O$.  Such a disk is a global attractor, since $\dot r < 0$ for $r > 1$.

\end{document}